# Infinite matrix of odd natural numbers. A bit about Sophie Germain prime numbers

Gennady Eremin

**Abstract.** In this paper we work with number sequences from the On-Line Encyclopedia of Integer Sequences (OEIS). Using the Pepis-Kalmar pairing function, we obtain an infinite matrix of natural numbers in which odd natural numbers are separated from even ones; such a matrix simplifies working with prime numbers. With point's shell numbers and shell lines we give another proof of the infinity of primes, based on Bertrand's postulate. We have shown that the asymptotic density of Mersenne numbers (OEIS A000225) is positive and equal to an infinitesimal value. It is also proven that Sophie Germain prime numbers (OEIS A005384) are located in the set of natural numbers with asymptotic density of $1/6$.

*Keywords*: OEIS, matrix of natural numbers, Mersenne number, pairing function, Bertrand's postulate, Ulam spiral, matrix carpet, Sophie Germain prime number.

## 1   Introduction

In this paper we work with some integer sequences from the well-known OEIS (the On-Line Encyclopedia of Integer Sequences [OEIS]). We will mainly be dealing with *primes*, i.e., natural numbers that are only multiples of 1 and the number itself (OEIS A000040). All primes are odd natural numbers greater than 1, except 2. In this regard, we will reform the set of natural numbers $\mathbb{N} = \{0, 1, 2, ...\}$ (like most mathematicians, we consider 0 to be a natural number, and this is also the ISO 80000-2:2019 standard). In the set $\mathbb{N}$, it is desirable to separate odd numbers from even ones, which will significantly simplify and speed up data processing when solving certain problems. For example, in various integer sequences it is often necessary to prove the infinity of prime numbers or primes of a special type, such as *Sophie Germain primes* (OEIS A005384), which we will also work with a bit in this paper. To transform the set of natural numbers, we use the *pairing function* $F : \mathbb{N}^2 \to \mathbb{N}$.

In the 19th century, mathematician Georg Cantor introduced the first pairing functions to uniquely encode two natural numbers into a single natural number, and these are two well-known functions of the following kind [Can78]:

(1) $\qquad \pi(x, y) = (x + y)(x + y + 1)/2 + y \quad$ and $\quad \pi'(x, y) = \pi(y, x).$

Cantor's goal was to prove that $\mathbb{N}^2$ is as countable as $\mathbb{N}$. It is proved that such quadratic functions are bijective and only polynomial pairing functions for the non-negative integers (see the Fueter-Pólya theorem in [Nat15, Vse02]).

Various variants of the Cantor pairing function are known in the literature, and in each variant we obtain some permutation of natural numbers in the first quadrant of the integer lattice (the right-top side region in the two dimensional space). However, in all such permutations, the even natural numbers are not separated from the odd natural numbers. In the middle of the 20th century, a pairing function appeared, which became popular and which gives us the desired permutation of natural numbers. The authors of this pairing function are the Polish mathematician Józef Pepis [Pep38] and the Hungarian mathematician Laszlo Kalmar [Kal39] (see [Tar08, Tar13]).

## 2 The Pepis-Kalmar pairing function

In [Ere24] we considered the Pepis-Kalmar pairing function in the following form:

(2) $$F(x, y) = 2^y (2x + 1) - 1 = z.$$

Function (2) is strictly monotone in each argument: for all $x, y \in \mathbb{N}$ we have both

(3)  a) $F(x, y) < F(x+1, y) = F(x, y) + 2^{y+1}$  and  b) $F(x, y) < F(x, y+1) = 2F(x, y) + 1$.

Also, function (2) is a bijection, since every positive integer $z+1$ can be written uniquely as the product of a power of two and an odd number.

Function (2) is of both practical and theoretical interest. Variations of this pairing function have been used by various authors in computer science and set theory, presumably because the inverse function $F^{-1}(z)$ is easily computed from the binary representation of $z+1$. In Figure 1, we have shown an initial fragment of the resulting matrix (prime numbers are shown in red).

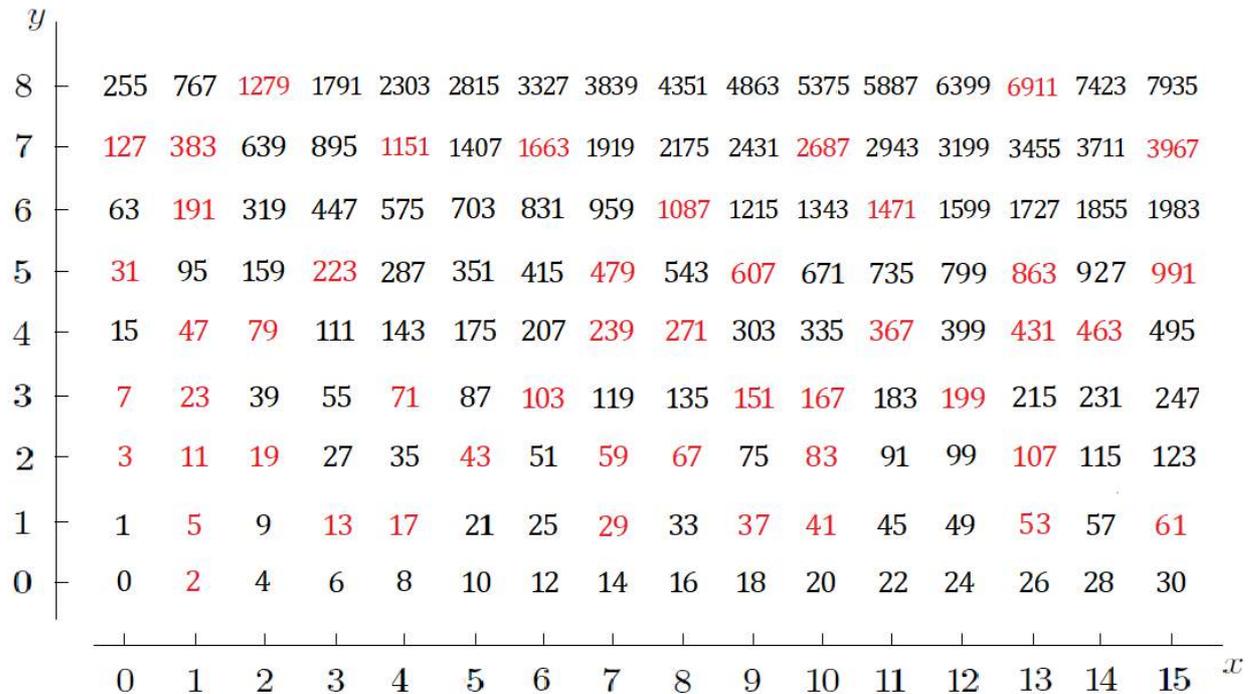

Figure 1. Infinite matrix of natural numbers (initial fragment).

As we can see, in the matrix even natural numbers are placed only on the x-axis, i.e., $F(x, 0) = 2x$ (OEIS A005843) or every even number $z = F(z/2, 0)$. In fact, it is the odd natural numbers that form a two-dimensional matrix, so such a construction can be called an *infinite matrix of odd natural numbers*. The y-axis contains Mersenne numbers, i.e., $F(0, y) = 2^y - 1 = M_y$ (*arithmetic-geometric progression* with the recurrent formula $a_{n+1} = 2a_n + 1$).

In the n-th row of the matrix (n-row), we obtain an *arithmetic progression* of the form $M_n + xd_n$ with the difference $d_n = 2^{n+1}$, and in each subsequent row the difference doubles (see (3a)). Arithmetic progressions do not intersect, and we will prove this below. The recurrence formula of Mersenne numbers (3b) works in the matrix columns, so we have arithmetic-geometric progressions in the columns. The columns also do not intersect and it is logical to call them *Mersenne columns* (or *unary Mersenne trees*). We will also prove the non-intersection of arithmetic-geometric progressions below. The first arithmetic progressions in rows, as well as the first arithmetic-geometric progressions in the matrix columns, we can find in the OEIS (see [Ere24]). The resulting matrix is convenient for working with odd prime numbers .

**2.1. Generating functions.** The generating function allows working with combinatorial objects by analytical methods. In mathematics, a *generating function* is an infinite formal power series with integer coefficients. It is believed that the sequence of such numbers is *generated* by a generating function. It is convenient when the generating function is obtained in closed form as of some expression.

In our matrix, the zero column $F(0, y) = 2^y - 1$ (the 0-column, the numbers with coordinate $x = 0$) contains the Mersenne numbers with generating function $B_0(y) = y/((1-y)(1-2y))$ (see OEIS A000225). The next 1-column contains the sequence of numbers $F(1, y) = 3 \times 2^y - 1$ with generating function $B_1(y) = (2-y)/((1-y)(1-2y))$ (OEIS A153893). Let us give still for the 2-column numbers $F(2, y) = 5 \times 2^y - 1$ the generating function $B_2(y) = (4-3y)/((1-y)(1-2y))$ (OEIS A153894). In general, for an arithmetic-geometric progression in an arbitrary $x$-column, we obtain the following generating function:

(4) $\qquad B_x(y) = (2x - (2x-1)y) / ((1-y)(1-2y)), \text{ for } x = 0, 1, 2, ...$

As a result, formula (4) gives us a family of generating functions that generate disjoint arithmetic-geometric progressions in the columns of the matrix. Let us recall that the union of such progressions is the set of all natural numbers.

Next, we will deal with arithmetic progressions in the rows of the matrix. On the $x$-axis we have even natural numbers, the sequence $F(x, 0) = 2x$ with the generating function $G_0(x) = 2x/(1-x)^2$ (OEIS A005843). In the next row ($y = 1$) we have an arithmetic progression $F(x, 1) = 1 + 4x$ with generating function $G_1(x) = (1+3x)/(1-x)^2$ (OEIS A016813). Still consider the arithmetic progression $F(x, 2) = 3 + 8x$ with the generating function $G_2(x) = (3+5x)/(1-x)^2$ (OEIS A017101), and the arithmetic progression $F(x, 3) = 7+16x$ with the generating function $G_3(x) = (7+9x)/(1-x)^2$ (there is no such sequence in the OEIS yet). As a result, we obtain another generating function for the Pepis-Kalmar pairing function:

(5) $\qquad G_y(x) = (M_y + (M_y+2)x) / (1-x)^2, \quad M_y = 2^y - 1, \text{ for } y = 0, 1, 2, \ldots$

Obviously, the family of formulas (5) generates an infinite family of generating functions for non-intersecting arithmetic progressions in the rows of the matrix, starting from the sequence of even natural numbers on the $x$-axis. The union of arithmetic progressions in the rows of the matrix as well as the union of arithmetic-geometric progressions (matrix columns) is the set of all natural numbers.

Thus, formula (4) and formula (5) generate the same set of integers (the natural series $\mathbb{N}$), but the author failed to find a connection between these two formulas.

**2.2. Asymptotic density**. In number theory, the *asymptotic density* (or *natural density*) means some value that allows us to estimate how large a particular set of natural numbers is. The density of the set of natural numbers $\mathbb{N}$ is 1, and accordingly the total asymptotic density of all rows in our matrix is also 1. The density of even natural numbers on the *x*-axis (0-row) is ½, since in $\mathbb{N}$ every second number is even. And then the density of odd natural numbers is 1 – ½ = ½ (every second natural number is odd); in our matrix this is the total density of all numbers in the *y*-rows, $y > 0$.

The row of odd numbers with coordinate $y = 1$ has an asymptotic density of ¼ (in the 1-row we have the arithmetic progression $4x +1$); then in each subsequent arithmetic progression the density decreases by 2 times, i.e., it is equal to half the density of the previous arithmetic progression. In general, according to formula (2), the density of the set of natural numbers with *y*-coordinate (*y*-row) is $2^{-y-1}$ (arithmetic progression $2^{y+1}x + M_y$). It is easy to see that in our matrix the asymptotic density of each arithmetic progression is equal to the sum of the densities of all subsequent arithmetic progressions (and obviously, each such sum is infinite). Thus, the following statement is true.

**Proposition 1.** *The asymptotic density of numbers in the y-row of the matrix is positive and equal to* $2^{-y-1}$.

It is easy to see that the asymptotic density of numbers in each column of the matrix approaches 0 as the corresponding interval of the numerical axis increases. It is enough to look up the values of the first 1000 Mersenne numbers (such data are in OEIS A000225): $M_{10} = 1023$ (the density is $10/1023 \sim 0.01$), $M_{50} = 1{,}125 \times 10^{15}$ (the density is $50/10^{15}$), $M_{500} = 3{,}373 \times 10^{150}$, $M_{1000} = 1{,}0715 \times 10^{301}$ (the density is $\sim 10^{-298}$), and so on.

As we can see, the number of Mersenne numbers in a finite interval of natural numbers rapidly decreases as the interval increases, i.e., in the long run we have an *infinitesimal positive density*. We observe the same in the other matrix columns, but at the same time we know that the total asymptotic density of arithmetic-geometric progressions in all columns is 1, since the columns (as well as the rows) contain all natural numbers. Hence it is easy to assume that in the matrix columns the asymptotic density of numbers is positive and equal to an infinitesimal value. In this case, the sum of infinitesimal values over all columns (of which there are an infinite number) gives us a density of 1.

In finite intervals, the density decreases as we move from column to column (as the *x*-coordinate increases). Therefore, for Mersenne numbers we are entitled to formulate the corresponding theorem, since in the initial column of the matrix the density is maximal in finite intervals.

**Theorem 2.** *The set of Mersenne numbers has non-zero natural density.*

The positive asymptotic density of Mersenne numbers may be significant, since we do not currently know whether there is a finite or infinite number of Mersenne primes (OEIS A000668). Obviously, the positive density of Mersenne numbers is necessary to prove the infinity of Mersenne primes (in fact, something similar happens in arithmetic progressions).

**2.3. Binary representation of integers.** We have considered two structures of the obtained natural matrix, these are arithmetic progressions in the rows of the matrix and arithmetic-geometric progressions in the columns of the matrix (we called the columns unary Mersenne trees, since the recurrent formula of Mersenne numbers works in each arithmetic-geometric progression). In each structure, the set of natural numbers $\mathbb{N}$ is partitioned into a family of infinite subsets (by rows or by columns).

In the resulting matrix, we often work with binary codes of numbers, for example, the recurrent formula of Mersenne numbers $a_{k+1} = 2a_k + 1$ corresponds to the simple operation of adding one unit to the end of the binary code of the original number $a_k$. Thus, the even number 6 (binary code 110) in the next row is transformed into $13 = 2 \times 6 + 1$ (code 1101); then in the following rows we get respectively the numbers 27 (11011), 55 (110111), and so on. As we can see, the numbers in the matrix columns are not repeated (as are the Mersenne numbers in the zero column). Since each column begins with a unique even natural number (even numbers on the $x$-axis are also not repeated), we can formulate the following obvious statement.

**Proposition 3.** *Each column of the matrix contains unique natural numbers (that is, the columns do not intersect).*

Next, let's work with strings, each of which starting with a unique Mersenne number. The difference in each row is more than twice the corresponding Mersenne number, so all numbers of such a row retain all the bits of the Mersenne number at the end of the binary code. For example, in a row that begins with the Mersenne number 15 (binary code 1111), the next number is $15 + 32 = 47$ (in binary code $1111 + 100000 = 101111$). At the second step we get the number $47 + 32 = 79$ (in binary code $101111 + 100000 = 1001111$), and further in all numbers of the row the binary suffix 01111 is preserved, since when adding the difference 32 only the prefix is increased. As a result, we get unique (non-repeating) numbers in the rows. Since each row starts with a unique Mersenne number (Mersenne numbers on the $y$-axis do not repeat), we can formulate another statement.

**Proposition 4.** *Each row of the matrix contains unique natural numbers (that is, the rows do not intersect).*

Clearly, Proposition 3 and Proposition 4 actually prove a bijection for the Pepis-Kalmar pairing function (2).

It is not hard to see, by the binary code of an arbitrary odd number $z = F(x, y)$ we simply determine the original coordinates $(x, y) = F^{-1}(z)$ (inverse pairing function): (a) the $y$-coordinate is equal to the number of 1's at the end of the binary code of the number $z$, and (b) we obtain the $x$-coordinate if we remove all final 1's and one 0 in the binary code of $z$ (if there are no zeros, then we have a Mersenne number, i.e., $x = 0$).

# 3    Shell-numbers and shell-lines in natural matrix

In this section, we consider another partitioning of the set of natural numbers, and this partitioning is directly related to the binary number codes. For pairing functions, mathematicians often introduce a point's *shell number*, which relates each point $(x, y)$ on an integer lattice to some value of $\sigma(x, y)$. For example, for the Cantor pairing function $\sigma(x, y) = x + y$, and then in the figure two neighboring points with the same shell number are connected by dots (or arrows, see [Szu19, Figure 1]). In this case the shell is a diagonal line, let us call it a *shell line* (or *shline*, for short); the set of points of such a line will be called a *shell set* (or *shset*, for short). As a result, all matrix points are partitioned into non-intersecting sets.

The American mathematician Matthew P. Szudzik proposed the following shell function for the Pepis-Kalmar pairing function [Szu19]:

(6) $$\sigma(x, y) = y + 1 + \lceil \log_2(x+1) \rceil,$$

where $\lceil t \rceil$ denotes the ceiling of $t$ for each real number $t$ (i.e., the nearest integer not less than $t$). In this case, according to equality (2) shell number $\sigma(x, y)$ is defined as the number of bits in the binary representation of $F(x, y) + 1$. Below Figure 2 shows the first five shlines (we borrowed Figure 3 from [Szu19]). The first shline (or 1-shline, for short) is point 0 with coordinates (0, 0), and it is the only point with the shell number 1. In all shell lines, the points with prime numbers are shown again in red.

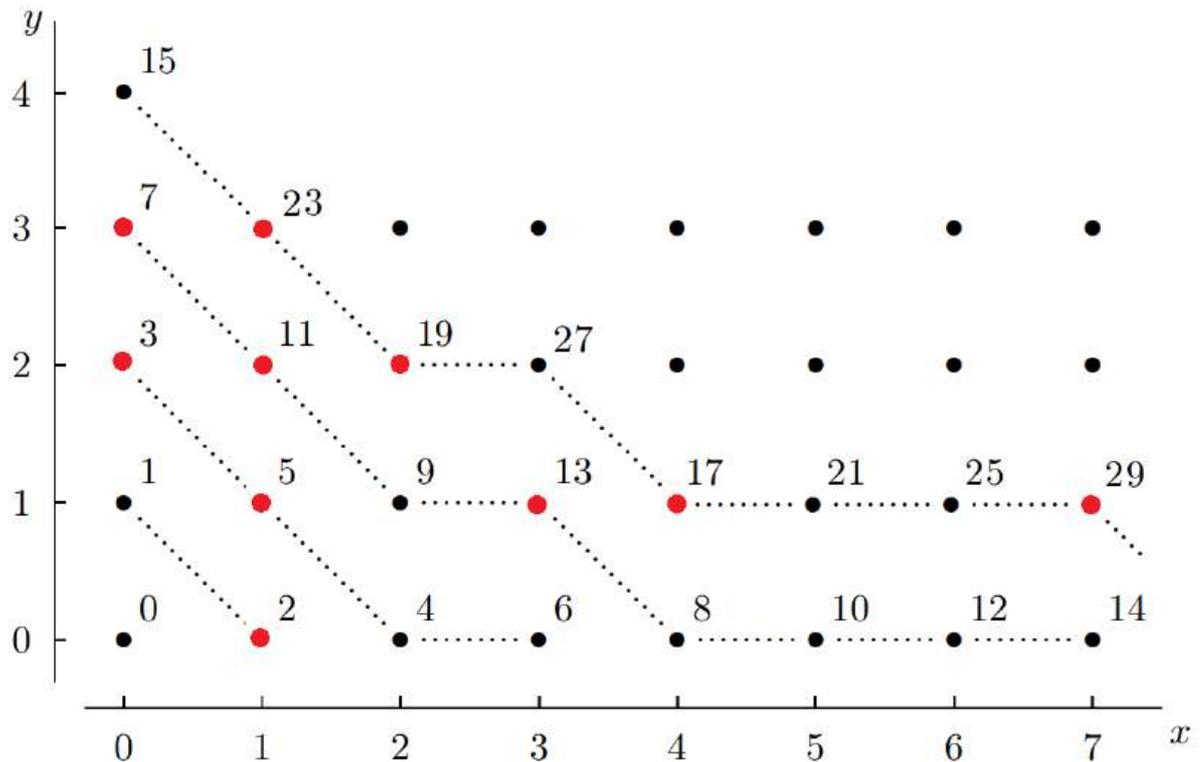

Figure 2. The points connected by the dotted line have the same shell number.

In each $(y+1)$-shline, the first term is the Mersenne number $M_y$, and it is the minimum number in such a shline; the next two terms are the numbers on the descending diagonal $F(1, y-1)$ and $F(2, y-2)$ for $y > 1$. The $(y+1)$-shline ends with a maximum number, and it is the doubled Mersenne number $2M_y$. The number of terms in the shline is even and equal to $2M_y - M_y + 1 = M_y + 1$. The other half of the numbers in the $(y+1)$-shline are even numbers arranged on the $x$-axis in ascending order: $M_y + 1$, $M_y + 3$, ..., $2M_y$.

In Appendix A, we give a Python-program that, given a Mersenne number, computes the sequence of terms of the corresponding shell set (6-shset is obtained in the program).

Thus, each $(y+1)$-shline contains all natural numbers from the closed interval $[M_y, 2M_y]$, and here we recall Bertrand's postulate (aka the Bertrand-Chebyshev theorem), which states that *for every natural number $n > 1$ there is a prime number between $n$ and $2n$*.

**3.1. Bertrand's postulate.** In our matrix, the $y$-axis contains the Mersenne numbers, and the number of these numbers is infinite. Each Mersenne number $M_y$ starts a $(y+1)$- shline, which includes all natural numbers from $M_y$ to $2M_y$. Therefore, according to Bertrand's postulate, there is at least one prime number in every $(y+1)$-shline, $y > 1$. In the matrix there are an infinite number of shlines and all shlines (and correspondingly all shsets) do not intersect, so we can formulate another proof of Euclid's theorem (the reader will find almost all proofs of Euclid's theorem in [Mes22]).

**Theorem 5.** *The set of prime numbers is infinite.*

Let's use Bertrand's postulate for our arithmetic progressions. In [Ere24], in the matrix rows we consider the *initial segments* (*short arithmetic progressions*) with terms from the OEIS A036991. In the *y*-row, the segment begins with the Mersenne number $M_y$, the number of terms in the segment is equal to $M_y + 1$, the last number of the segment is composite and is equal to $M_y(2M_y + 3) = M_y(M_{y+1} + 2)$. Thus, the initial segment and the corresponding shline start with the same Mersenne number and have the same length.

Additionally note, in each shline the second half of the numbers are even, but in short arithmetic progressions there are no even numbers at all, so in each segment the area for searching for prime numbers is twice as large as in the corresponding shline. Given that the last number in the initial segments is composite (we reduce the interval for finding prime numbers by 1), we can formulate a hypothesis similar to Bertrand's postulate.

**Conjecture 6.** *In the y-row of the matrix, $y > 1$, in the short arithmetic progression that begins with the Mersenne number $M_y$ and has length equal to $M_y$, there is a prime number.*

Apparently, when proving Conjecture 6, we will have to consider the difference of the arithmetic progression in the *y*-row $d_y = (M_y(2M_y + 3) - M_y)/M_y = 2M_y + 2 = 2^{y+1}$.

Obviously, in the case of the proof of Conjecture 6 we automatically obtain the *infinity of prime numbers* in the OEIS A036991. Additionally, in the matrix rows we obtain the reduction of the known Linnik's constant to 2 for arithmetic progressions (see [Ere24]).

**3.2. Bijection.** The length of each shell line is equal to the length of the corresponding initial segment consisting of A036991 terms, so we get a one-to-one correspondence between these finite sets. In other words, we obtain a bijection between the set of natural numbers (including 0) and some subset of the OEIS A036991. Let us formulate the corresponding statement.

**Proposition 7.** *There is a bijection between the set of natural numbers and the short arithmetic progressions (initial segments) in the matrix rows.*

Below we give a small table of such a bijection for initial numbers. The first row of the table shows the terms of the first shlines (in order, starting with the Mersenne numbers), and this row contains all natural numbers, including 0. The second row of the table shows the terms of the corresponding short arithmetic progressions, and these are the terms of a subset of the sequence A036991 (recall that the sequence A036991 is also a subset of the natural numbers).

| 0 | 1 | 2 | 3 | 5 | 4 | 6 | 7 | 11 | 9 | 13 | 8 | 10 | 12 | 14 | 15 | 23 | 19 | 27 | ... |
|---|---|---|---|---|---|---|---|----|---|----|---|----|----|----|----|----|----|----|-----|
| 0 | 1 | 5 | 3 | 11 | 19 | 27 | 7 | 23 | 39 | 55 | 71 | 87 | 103 | 119 | 15 | 47 | 79 | 111 | ... |

For example, the second row does not contain the following A036991 terms: 13, 21, 29, 43, 45, 51, 53, 59, 61, 75, 77, ... ; also, the second row does not contain numbers that are not A036991 terms: 9, 17, 25, 33, 35, 37, 41, 49, 57, 65, 67, 69, ... Let's add, the second row has no even natural numbers at all (except 0).

In such a table we can consider the Mersenne numbers as control points. If we number the columns of the table (starting at 0), the value of each Mersenne number (as well as the value of the preceding number in the first row) will coincide with the ordinal number of the column. In Appendix B, we give a table for the first 127 columns.

## 4   Matrix carpet. A bit on Sophie Germain prime numbers

Often in drawings prime numbers are highlighted in a special way (color or otherwise), up to the elimination of composite numbers, in order to get a clear picture of the structure of prime numbers. For example, this is done in *Ulam spiral* (or *Ulam carpet*) [Ulam]. Below in Figure 3 we have selected the prime numbers of the matrix and shown them as dots, as it is usually done in the Ulam spiral. Let's call this pattern the *matrix carpet*.

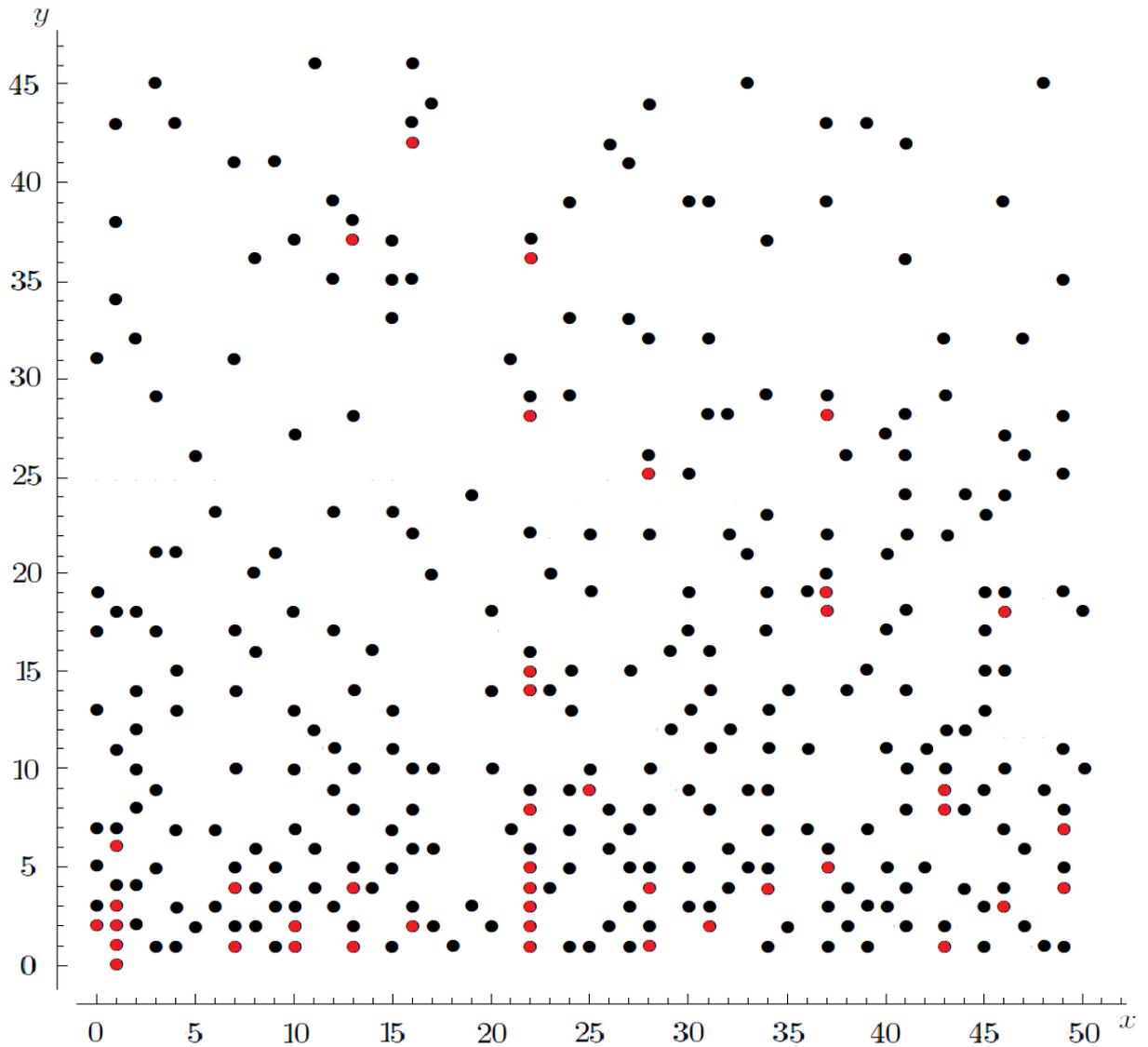

Figure 3. Matrix carpet, Sophie Germain primes shown in red.

In mathematics, there are many papers that study *left shifted primes* and *right shifted primes*. In our matrix, this is also available, but there is an additional feature, here it is convenient to work with prime numbers that are s*hifted up* (i.e., *up shifted primes*). Such are the *Sophie Germain prime numbers* and the associated s*afe prime numbers*. Recall that Sophie Germain prime is a prime number $p$ such that the number $2p+1$ is also prime; in this case, the number $2p+1$ is called a safe prime number.

In Figure 3 we have highlighted the Sophie Germain prime numbers in red, the first such numbers are 2, 3, 5, 11, ... (see OEIS A005384). As we see among the Mersenne numbers (0-column of the matrix) there is only one Sophie Germain prime number, it is 3 and the corresponding safe number is 7. And it is easy to explain, because every second Mersenne number after 3 is composite. This is proved simply: for $q = 2p+1$ the next number $2q+1 = 4p+3$ is a multiple of 3 if $p$ is a multiple of 3. Thus in the sequence of Mersenne numbers there are no two neighboring prime numbers greater than 3.

Obviously, in our case, according to formula (2), a prime number $F(x, y)$ is called Sophie Germain prime number if the number $F(x, y+1)$ is also prime. And then the second number is a safe prime number.

Note, in Figure 3, there are no Sophie Germain primes in many columns, and this is also not difficult to explain. The 1-column (points with coordinate $x = 1$) starts with the Sophie Germain prime numbers 2, 5, 11, 23 plus the safe prime number 47 (see also Figure 2). In the next 2-column we have a sequence of points 4, 9, 19, 39, 79, .... where every second number is a multiple of 3 (as in the case of Mersenne numbers). And this is because the initial number $4 \equiv 1 \pmod 3$, and the subsequent number $4 \times 2 + 1 = 9 \equiv 0 \pmod 3$, i.e. a multiple of 3. Next $9 \times 2 + 1 = 19 \equiv 1 \pmod 3$, and we again get 1 modulo 3. Thus, in such a sequence every second number is a multiple of 3, so two neighboring prime numbers are impossible in the 2-column.

Further, in the 3-column with the numbers 6, 13, 27, 55, ... there are no Sophie Germain prime numbers either, because such column starts with number 6, which is divisible by 3, and therefore every second number is also divisible by 3 (as in the previous 2-column). As a result, we can formulate another statement.

**Proposition 8.** *Sophie Germain prime numbers are possible only in those columns of the matrix where each number gives a remainder of 2 modulo 3.*

*Proof* is obvious, since if any natural number $p \equiv 2 \pmod 3$, then also the subsequent number $q = 2p+1 \equiv 2 \pmod 3$. That is, in such columns there are no numbers divisible by 3 at all, and so neighboring prime numbers are possible. □

As a result, we conclude that Sophie Germain prime numbers are possible only in the columns of the matrix with coordinate $x = 1, 4, 7, 10, ...$, and these columns begin with even numbers respectively 2, 8, 14, 20, ..., i.e., numbers of the form $6k + 2$, each of which gives a remainder 2 modulo 3 (OEIS has the corresponding sequence A016933). The remaining columns of the matrix can be considered as a forbidden zone for searching for Sophie Germain prime numbers.

Thus, Sophie Germain prime numbers may appear in only one-third of the columns of the matrix. Since we are working with odd natural numbers (the asymptotic density is ½), the search for Sophie Germain prime numbers is effectively limited to the set of natural numbers with asymptotic density (½)/3 = ¹/₆. And we can formulate the following theo-rem.

**Theorem 9**. *Sophie Germain prime numbers are placed in the set of natural numbers with asymptotic density equal to $1/6$.*

In mathematics, it is assumed that the number of Sophie Germain prime numbers is infinite, but this is an open question in number theory. It is possible that narrowing the search area for Sophie Germain prime numbers will make it easier to prove their infinity in the future.

## Acknowledgements.

The author is very grateful to Olena Kachko and Igor Shparlinski for many discussion and useful comments on preliminary versions of the paper. Also the author would like to thank the anonymous referees for their valuable comments and suggestions.

[Vse02]   M. A. Vsemirnov. *Two elementary proofs of the Fueter-Pólya theorem on pairing polynomials*, St. Petersburg Mathematical Journal, **13**(5), 705-715 (2002).

[OEIS]   The On-Line Encyclopedia of Integer Sequences, 2024.   https://oeis.org/

Concerned with sequences A000040, A000225, A000668, A005384, A005843, A016813, A017101,  A036991, A153893,  A153894.

*Email address*:  ergenns@gmail.com

Written:  January 28, 2025

## Appendix A.   Python-program for calculating a shell set

```
#  we get a shell set by a given shell number
shnumb = 6     # better  shell number > 3
mers= 2**(shnumb - 1) - 1; shset= [mers]
diff = 2*(mers + 1); term = mers
while 1:
   if term + diff < 2*mers: term += diff
   else:
      diff //= 2; term = (term - 1)//2 + diff
   if diff == 2: break
   shset.append(term)
print(shset)
print(f'Next come even numbers: {term}, {term+2}, ..., {2*mers}. ')
```

The result is as follows:  31, 47, 39, 55, 35, 43, 51, 59, 33, 37, 41, 45, 49, 53, 57, 61.
Next come even numbers:  32, 34, 36, ...,  62.

# Appendix B. Bijection: shell lines and short arithmetic progressions

| 0 | 1 | 2 | 3 | 5 | 4 | 6 | 7 | 11 | 9 | 13 | 8 | 10 | 12 | 14 | 15 | 23 | 19 | 27 | 17 |
|---|---|---|---|---|---|---|---|---|---|---|---|---|---|---|---|---|---|---|---|
| 0 | 1 | 5 | 3 | 11 | 19 | 27 | 7 | 23 | 39 | 55 | 71 | 87 | 103 | 119 | 15 | 47 | 79 | 111 | 143 |

| 21 | 25 | 29 | 16 | 18 | 20 | 22 | 24 | 26 | 28 | 30 | 31 | 47 | 39 | 55 |
|---|---|---|---|---|---|---|---|---|---|---|---|---|---|---|
| 175 | 207 | 239 | 271 | 303 | 335 | 367 | 399 | 431 | 463 | 495 | 31 | 95 | 159 | 223 |

| 35 | 43 | 51 | 59 | 33 | 37 | 41 | 45 | 49 | 53 | 57 | 61 | 32 | 34 |
|---|---|---|---|---|---|---|---|---|---|---|---|---|---|
| 287 | 351 | 415 | 479 | 543 | 607 | 671 | 735 | 799 | 863 | 927 | 991 | 1055 | 1119 |

| 36 | 38 | 40 | 42 | 44 | 46 | 48 | 50 | 52 | 54 | 56 | 58 |
|---|---|---|---|---|---|---|---|---|---|---|---|
| 1183 | 1247 | 1311 | 1375 | 1439 | 1503 | 1567 | 1631 | 1695 | 1759 | 1823 | 1887 |

| 60 | 62 | 63 | 95 | 79 | 111 | 71 | 87 | 103 | 119 | 67 | 75 | 83 | 91 |
|---|---|---|---|---|---|---|---|---|---|---|---|---|---|
| 1951 | 2015 | 63 | 191 | 319 | 447 | 575 | 703 | 831 | 959 | 1087 | 1215 | 1343 | 1471 |

| 99 | 107 | 115 | 123 | 65 | 69 | 73 | 77 | 81 | 85 | 89 | 93 | 97 |
|---|---|---|---|---|---|---|---|---|---|---|---|---|
| 1599 | 1727 | 1855 | 1983 | 2111 | 2239 | 2367 | 2495 | 2623 | 2751 | 2879 | 3007 | 3135 |

| 101 | 105 | 109 | 113 | 117 | 121 | 125 | 64 | 66 | 68 | 70 | 72 | 74 |
|---|---|---|---|---|---|---|---|---|---|---|---|---|
| 3263 | 3391 | 3519 | 3647 | 3775 | 3903 | 4031 | 4159 | 4287 | 4415 | 4543 | 4671 | 4799 |

| 76 | 78 | 80 | 82 | 84 | 86 | 88 | 90 | 92 | 94 | 96 | 98 | 100 |
|---|---|---|---|---|---|---|---|---|---|---|---|---|
| 4927 | 5055 | 5183 | 5311 | 5439 | 5567 | 5695 | 5823 | 5951 | 6079 | 6207 | 6335 | 6463 |

| 102 | 104 | 106 | 108 | 110 | 112 | 114 | 116 | 118 | 120 | 122 | 124 | 126 |
|---|---|---|---|---|---|---|---|---|---|---|---|---|
| 6591 | 6719 | 6847 | 6975 | 7103 | 7231 | 7359 | 7487 | 7615 | 7743 | 7871 | 7999 | 8127 |